\documentclass[12pt]{article}
\usepackage{amsmath,amsfonts,amssymb,amscd}
\usepackage{graphicx}
\usepackage[T2A]{fontenc}

\newtheorem{theo}{Theorem}
\newtheorem{lem}{Lemma}[section]
\newtheorem{prop}{Proposition}[section]
\newtheorem{conj}{Conjecture}
\newtheorem{cor}{Corollary}[section]

\newtheorem{dfn}{Definition}[section]

\makeatletter \@addtoreset{equation}{section} \makeatother

\newcommand{\pr}{\mathop{\rm pr}\nolimits}

\newcommand{\mC}{\mathbb{C}}

\newcommand{\mR}{\mathbb{R}}
\newcommand{\mT}{\mathbb{T}}
\newcommand{\mZ}{\mathbb{Z}}

\newcommand{\bB}{{\bf B}}

\newcommand{\bI}{{\bf I}}

\newcommand{\ba}{{\bf a}}
\newcommand{\bb}{{\bf b}}
\newcommand{\bc}{{\bf c}}

\newcommand{\bn}{{\bf n}}

\newcommand{\DD}{{\cal D}}
\newcommand{\FF}{{\cal F}}

\newcommand{\PP}{{\cal P}}

\newcommand{\TT}{{\cal T}}

\newcommand{\thet}{\vartheta}

\newcommand{\meas}{\operatorname{meas}}
\newcommand{\id}{\operatorname{id}}

\newcommand\qed{{\unskip\nobreak\hfil\penalty50
  \hskip2em\hbox{}\nobreak\hfil\mbox{\rule{1ex}{1ex} \qquad}
    \parfillskip=0pt \finalhyphendemerits=0\par\medskip}}

\begin{document}

\title
{A locally integrable multi-dimensional billiard system}
\author{D.Treschev}
\date{}
\maketitle


\section{Introduction}
\label{sec:intro}

Let $D\subset E$ be a domain with the smooth boundary $S = \partial D$ in the Euclidean space $E = \mR^{n+1}$, $n\ge 1$. Assuming that the closure $\bar D$ is compact, we define the billiard system in $D$ as follows. A particle moves with a unit constant velocity inside $D$. The reflection from the boundary is elastic.

The word elastic means a well-known relation between the velocities of the particle before and after the impact (``the angle of the incidence equals the angle of the reflection''), but we will use an equivalent variational equation.
Consider the generating function (the discrete Lagrangian)
$$
  L:\PP\to\mR, \quad
  \PP = S\times S, \qquad
  L(\ba,\bb) = |\ba-\bb|.
$$
Then $\ba,\bb,\bc\in S$ are 3 consecutive impact points iff
\begin{equation}
\label{Fermat}
    \partial_\bb \big( L(\ba,\bb) + L(\bb,\bc) \big)
  = 0, \qquad
    \partial_\bb
  = \frac{\partial}{\partial\bb}.  
\end{equation}

Billiard systems were introduced by Birkhoff \cite{Birkhoff} and since that time occupied a noticeable part of the dynamics. Many results and references on the subject, as a rule in the case $n=1$, can be found in \cite{KT,Tab}.

Now we define the billiard map $\beta:\PP\to\PP$. A pair of consecutive impact points $(\ba,\bb)$ is transformed to
$(\bb,\bc) = \beta(\ba,\bb)$ such that $\ba,\bb,\bc$ satisfy (\ref{Fermat}).
This means that $\beta$ determines a discrete Lagrangian system (a general definition of a discrete Lagrangian system is contained in \cite{BT_anti}). In particular, $\beta$ is symplectic.

In this paper we consider the case when the domain $D$ is symmetric in all coordinate hyperplanes. We assume that $E$ splits into the direct product $\mR^n\times\mR$, where the subspaces $\mR^n\times\{0\}\subset E$ and $\{0\}\times\mR\subset E$ are called horizontal and vertical respectively. The map $\beta$ has the periodic orbit $\gamma$ of period 2 lying on the vertical coordinate axis. We are interested in the local dynamics near $\gamma$.

Consider the ball $B = \{x\in\mR^n : |x| < d\}$. We determine $S$ locally by the graphs
$$
  S_- = \{(x,f(x)) : x\in B \}, \qquad
  S_+ = \{(x,-f(x)) :  x\in B \},
$$
where
\begin{equation}
\label{f}
 \mbox{$f : B\to\mR$
       is a negative even function.}
\end{equation}
Our question is as follows.
\medskip

{\bf Q}. Is it possible to choose $f$ so that the corresponding billiard map $\beta$ is locally (near $\gamma$) conjugated to a linear map?
\medskip

Before a discussion of the motivations and results, we reformulate this question in a technically more convenient form.

Let $I: E\to E$, $I^2 = \id$ be the symmetry in the horizontal hyperplane. Then $I$ can be naturally extended to an involution of the phase space $\PP$: $(\ba,\bb)\mapsto (I\ba,I\bb)$. Slightly abusing the notation, we denote this involution also as $I$. We define the quotient space $\hat\PP = \PP / I$ and the natural projection $\pr:\PP\to\hat\PP$.

The maps $I,\beta:\PP\to\PP$ commute. Therefore there exists a unique map
$\hat\beta : \hat\PP\to\hat\PP$ such that
$$
  \hat\beta\circ\pr = \beta.
$$
The quotient map $\hat\beta$ is more convenient for the local study because $\gamma$ projects to a fixed point $\hat\gamma$ of $\hat\beta$. Hence we can propose another (equivalent) version of the question {\bf Q}.
\medskip

{$\bf\hat Q$}. Is it possible to choose $f$ so that the corresponding quotient billiard map $\hat\beta$ is locally (near $\hat\gamma$) conjugated to a linear map?
\medskip

This linear map $\rho$ is symplectic and coincides with the linearization of $\hat\beta$ at $\hat\gamma$. Let $\lambda_1^{\pm 1},\ldots,\lambda_n^{\pm 1}$ be the eigenvalues of $\rho$.

In this paper we consider only the most interesting situation when $\hat\gamma$ is linearly stable. Then the eigenvalues lie on the unit circle. Explicit conditions of linear stability of a billiard orbit of period two (in general, not necessarily symmetric case) are well-known in the case $n=1$ (see for example, \cite{KT}), some geometric interpretations of these conditions are given in \cite{K2d}. In the case $n=2$ such stability conditions can be found in \cite{K3d}.

If the eigenvalues form a resonant set i.e.,
$$
  \lambda_1^{k_1} \cdots \lambda_n^{k_n} = 1 \quad
  \mbox{for some nonzero vector 
         $k=(k_1,\ldots,k_n)\in\mZ^n$},
$$
then there is no hope to have a positive answer to ${\bf\hat Q}$ (see Section \ref{sec:formal}).  In the nonresonant case we show (Theorem \ref{theo:formal}) that $f$ can be obtained as a formal Taylor series. The most intriguing question on the convergence of this series remains open. But numerical analysis makes reasonable the following

\begin{conj}
\label{conj:main}
If the set $\lambda_1,\ldots,\lambda_n$ is nonresonant and moreover, satisfies good Diophantine properties, the series, presenting $f$, is locally convergent. The same is true for the series, presenting the conjugacy map.
\end{conj}

Consider the nonresonant case. If Conjecture \ref{conj:main} is true, we obtain a local real-analytic billiard system with linear quasi-periodic dynamics. In this case $f$ turns out to be convex. This locally defined convex function $f$ can be smoothly continued up to a function $F$ such that graphs of $\pm F$ determine a smooth closed hypersurface $S\subset E$, the boundary of a convex domain $D\subset E$. Hence we obtain a globally defined billiard system such that the corresponding phase space $\PP$ contains an open domain $\DD\subset\PP$ filled by quasiperiodic motions with the same set of frequencies. In particular, the only periodic orbit in $\DD$ is $\gamma$.

Billiard system in the domain $D$ is locally integrable. The problem of integrability of billiard systems is widely discussed. According to the Birkhoff conjecture any domain on a plane bounded by a smooth closed curve generates an integrable billiard system iff this curve is an ellipse. Partial results confirming this conjecture are contained in \cite{Bol_bill,Bia-Mir,ASK}, (see also \cite{Glu_bill}, where an analog of the Birkhoff conjecture for outer billiard systems is proven). Usually the billiard integrability is discussed in the context of the existence of a first integral polynomial in the momenta, see a recent survey and a collection of new results in \cite{K_bill_UMN}. Local first integrals which should exist for the billiard system in $D$ are real-analytic, but certainly not polynomial in the momenta.

Because of a strong degeneracy of the billiard dynamics in $D$ we expect that the spectrum of the Laplace operator in $D$ (say, with the Dirichlet boundary conditions) can be very special. Thereby an interesting question appears on possible values of the quantity $\meas\DD / \meas\PP$. Numeric computations in the case $n=1$ show that this ratio can exceed 50\%.

Now we present some discussion of numeric results for $n=1$ and $n=2$.
\medskip

{\bf The case $n=1$}. This situation is studied in \cite{Tre_PhysD} and \cite{Tre_Proc}. We consider the normalization $f(0)=-1/2$. Taking $\lambda=\lambda_1 = e^{i\alpha}$, where $\alpha/\pi$ is an irrational number, we compute the sequence of coefficients $f_{2j}$, where
\begin{equation}
\label{f=sum}
  f(x) = \sum f_{2j}x^{2j}.
\end{equation}
Now we present some conjectures motivated by results of numeric computations.

1. For Diophantine $\alpha/(2\pi)\in (0.3,0.5)$ the ratio
   $b_j = f_{2j} / f_{2j-2}$ admits the asymptotic expansion
$$
  b_j = b_\infty\Big(1 + \frac{\sigma}{j} + O\Big(\frac{1}{j^2}\Big)\Big).
$$

\begin{figure}[h]
\begin{center}
\vskip-2mm
\includegraphics[scale=0.5]{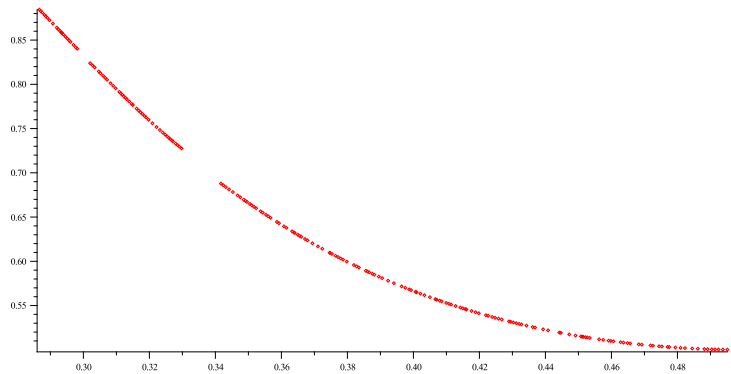}
\caption{The graph of $b^{-1/2}_\infty$ as a function of $\alpha/(2\pi)$. Two ``gaps'' correspond to the resonances $\frac\alpha{2\pi} = 3/10$ and
$\frac\alpha{2\pi} = 1/3$.}
\label{fig:binfty}
\end{center}
\end{figure}

2. Graph of $b_\infty^{-1/2}$ as a function of $\alpha/(2\pi)$ is presented in Fig. \ref{fig:binfty}.

This function is not defined for rational $\alpha/\pi$, but looks smooth. Probably, this function is Whitney smooth, \cite{Whitney}.

3. Independently of $\alpha$ numerically $\sigma = -3/2$. If this really the case, we have:
$$
  f_{2j} = C \frac{b_\infty^j}{j^{3/2}} \Big(1 + O\Big(\frac1j\Big)\Big).
$$
This would imply that series (\ref{f=sum}) converges on the boundary of the convergence disk $|x|\le x_* = b_\infty^{-1/2}$ and has at the points
$x = \pm x_*$ singularities of type $\sqrt{(x_* \mp x)}$, in particular, the tangent line to the graph at these points is vertical.

This means that the graph of $f$ can be probably continued through the points $(\pm x_*,f(\pm x_*)$ up to a longer real-analytic curve.
\medskip

{\bf The case $n=2$}. Putting $f(0) = -1/2$, we compute the Taylor coefficients $a_{j_1 j_2}$ ($j_1,j_2$ are even), where
$$
    f(x)
  = \frac12 \sum_{j_1,j_2\in 2\mZ_+} a_{j_1 j_2} x_1^{j_1} x_2^{j_2}, \qquad
    a_{00} = -1, \quad
    \mZ_+ = \{0,1,\ldots\} .
$$

The coefficients $a_{0k}$ and $a_{k0}$ can be computed from the case $n=1$ because sections of the billiard domain by the vertical planes ${x_1=0}$ and ${x_2=0}$ give solutions of the problem with $n=1$.

Take for example, $\lambda_1 = e^{i\alpha_1}$, $\lambda_2 = e^{i\alpha_2}$, where $\alpha_1/(2\pi)$ and $\alpha_2/(2\pi)$ are quadratic irrationals:
$$
  \alpha_1 = 2\pi\cdot(3,3,1,1,1,1,\ldots),\quad
  \alpha_2 = 2\pi\cdot(2,5,2,2,2,2,\ldots)
$$
(chain fractions). For any even $k$ we present the line
$$
  \frac{a_{0,k}}{\sqrt{C_k^0}},\frac{a_{2,k-2}}{\sqrt{C_k^2}}, \ldots ,
  \frac{a_{k,0}}{\sqrt{C_k^k}}.
$$
The binomial multipliers $1/\sqrt{C_k^l}$ are motivated by the Bombieri norm on the space of homogeneous polynomials \cite{Bomb}. Here are numeric data beginning from $k=4$ (we save only 5 digits):
\begin{eqnarray*}
& .50276,\;\, 1.0749,\;\, 1.8853 &\\
& .38788,\;\, 1.1811,\;\, 1.9557,\;\, 3.6123 &\\
& .36853,\;\, 1.5808,\;\, 2.7866,\;\, 4.5700,\;\, 8.6479 &\\
& .39228,\;\, 2.3233,\;\, 4.4113,\;\, 7.3709,\;\, 12.080,\;\, 23.183&\\
& .44643,\;\, 3.6066,\;\, 7.3683,\;\, 12.798,\;\, 20.965,\;\, 34.380,
     \;\, 66.587 &\\
& .53202,\;\, 5.8039,\;\, 12.711,\;\, 23.049,\;\, 38.630,\;\, 62.628,
     \;\, 102.77,\;\, 200.34
\end{eqnarray*}
Conjecture \ref{conj:main} for $n=1$ predicts an exponential growth of the numbers $a_{0k}$ (the left-hand side of the table) and $a_{k0}$ (the right-hand side of the table). We see that the numbers in each line grow monotonically. This suggests an extension of Conjecture \ref{conj:main} to the case of arbitrary $n\ge 1$.

The further plan of the paper is as follows. In Section \ref{sec:conj} we obtain equations (the conjugacy equations) from which the function $f$ and the conjugacy map can be computed. We discuss basic symmetries of this equation and some properties of its (formal) solutions in Section \ref{sec:symm}. A further analysis of the solutions is contained in Section \ref{sec:formal}, where we prove a theorem about the existence of a formal solution in the non-resonant case. Other symmetries of the formal solution are discussed in Section \ref{sec:other}. Finally we present another form of the conjugacy equation which contains the unknown functions polynomially. We used this equation in the numeric analysis of the Taylor series for $f$ and the conjugacy map.

\section{Conjugacy equation}
\label{sec:conj}

\subsection{Conjugacy map}

We consider billiard trajectories which hit $S_-$ and $S_+$ alternatively. If $\ba$ and $\bb$ are two consequtive impact points of a trajectory then
$$
  \ba = \Big(\begin{array}{c} a\\ f(a)\end{array}\Big), \;
  \bb = \Big(\begin{array}{c} b\\ -f(b)\end{array}\Big) \quad
  \mbox{or} \quad
  \ba = \Big(\begin{array}{c} a\\ -f(a)\end{array}\Big), \;
  \bb = \Big(\begin{array}{c} b\\ f(b)\end{array}\Big), \qquad
  a,b\in B.
$$
In both cases
\begin{equation}
\label{discr_L}
    L(\ba,\bb)
  = |\ba - \bb|
  = \hat L(a,b)
  = \Big( |a-b|^2 + (f(a) + f(b))^2 \Big)^{1/2}.
\end{equation}
Here the passage from the points $\ba,\bb\in S$ to their projections to $B$ corresponds to the passage from $\beta$ to the quotient map $\hat\beta$. Note that in these coordinates 
$$
  \hat\gamma = (0,0)\in B\times B.
$$

Hence $\ba,\bb,\bc$ are 3 consequtive points on a trajectory\footnote
{which hits $S_\pm$ alternatively}
of $\beta$ iff $a,b,c$ are 3 consequtive points on a trajectory of $\hat\beta$ iff (\ref{Fermat}) holds or equivalently,
$$
  \frac{\partial}{\partial b}
    \big( \hat L(a,b) + \hat L(b,c) \big) = 0. 
$$
More explicitly,
\begin{equation}
\label{Fermat_expl}
    \frac{b-a + \nabla f(b)(f(a)+f(b))}{L(a,b)}
  + \frac{c-a + \nabla f(b)(f(c)+f(b))}{L(b,c)}
  = 0.
\end{equation}

Given a collection of complex numbers
\begin{equation}
\label{lambda}
  \lambda_1,\ldots,\lambda_n, \qquad
  |\lambda_1| = \ldots = |\lambda_n| = 1
\end{equation}
consider the linear symplectic map
$$
  \rho:\mC^n\to\mC^n, \quad
          z = (z_1,\ldots,z_n)
  \mapsto \rho(z) = (\lambda_1 z_1,\ldots,\lambda_1 z_n).
$$

The billiard map is locally conjugated with $\rho$ iff there exists a diffeomorphism $X:U\to\bB$, where $U\subset\mC^n$ and $\bB\subset B\times B$ are neighborhoods of zero,\footnote
{Below we assume that $U$ is $\rho$-invariant.}
such that the following diagram
\begin{equation*}
\begin{CD}
     U   @>\rho>>   U \\
    @V{X}VV       @VV{X}V\\
    \bB  @> \hat\beta >>     \bB
\end{CD}
\end{equation*}
commutes. The map $X$ has no relation with the complex structure on $\mC^n$. Hence we have to use $(z,\bar z) = (z_1,\ldots,z_n,\bar z_1,\ldots,\bar z_n)$ as coordinates on $\mC^n$. We put
$$
  X(z,\bar z) = (a,b) = (\chi_-(z,\bar z),\chi(z,\bar z)).
$$

Let $\beta(\ba,\bb) = (\bb,\bc)$ and $a = \chi_-(z,\bar z)$. Then
$$
  b = \chi(z,\bar z) = \chi_-\circ\rho(z,\bar z), \quad
  c = \chi\circ\rho(z,\bar z),
$$
and $a,b,c$ are connected by (\ref{Fermat_expl}). This implies
\begin{equation}
\label{chichi}
  \chi = \chi_-\circ\rho
\end{equation}
and
\begin{equation}
\label{conj_Lagr}
     \partial_2 L(\chi\circ\rho^{-1},\chi)
   + \partial_1 L(\chi,\chi\circ\rho)
  = 0,
\end{equation}
where $\partial_1$ and $\partial_2$ are defined by
$$
  \partial_2 L(x,y) = \frac{\partial L}{\partial y}(x,y), \quad
  \partial_1 L(x,y) = \frac{\partial L}{\partial x}(x,y).
$$
In principle, the functions $f,\chi$ can be computed from (\ref{conj_Lagr}). But (\ref{conj_Lagr}) turns out to be overdetermined w.r.t. $f$. This produces some difficulties in the proof of the existence of a solution. To avoid these difficulties, below we replace (\ref{conj_Lagr}) by an equivalent system, free of such problems.

\subsection{Averaging}

We define the action of the torus $\mT^n = \mR^n / (2\pi\mZ^n)$ on $\mC^n$ by the equation
$$
  \mT\ni\alpha\mapsto \rho_\alpha : \mC^n \to\mC^n, \qquad
  \rho_\alpha z = e^{i\alpha} z.
$$

For any function $g:U\to\mC$ defined on a set $U\subset\mC^n$, invariant with respect to the action of $\mT^n$, we put $\rho^*_\alpha g = g\circ\rho_\alpha$ and define the average
$$
  \langle g\rangle:U\to\mC, \qquad
     \langle g\rangle
  = \frac1{(2\pi)^n}\int_{\mT^n} \rho^*_\alpha g\, d\alpha, \quad
    [g] = g - \langle g\rangle.
$$
For any $g = \sum_{k',k''\in\mZ_+^n} g_{k'k''} z^{k'}\bar z^{k''}$ we have the identities
$$
  \langle g\rangle = \sum_{k\in\mZ_+^n} g_{k k} (z\bar z)^k, \quad
  \langle g\rangle = \langle \rho^* g\rangle.
$$
Averaging can be applied to differential forms as well. In particular, we have the identity
\begin{equation}
\label{<>=<>}
  d\langle g\rangle = \langle dg\rangle \quad
  \mbox{for any function $g:U\to\mC$}.
\end{equation}

\subsection{The forms $\nu,\hat\nu$}

Consider the forms $\mu,\hat\mu$ on $B\times B$:
$$
    \mu
  = \partial_2 L(x_-,x)\, dx
  = \sum_{j=1}^n \partial_{x_j} L(x_-,x)\, dx_j ,  \quad
    \hat\mu
  = \partial_1 L(x_-,x)\, dx_-
  = \sum_{j=1}^n \partial_{x_{j-}} L(x_-,x)\, dx_{j-}.
$$
Then $d\mu = -d\hat\mu$ is the standard symplectic structure for the billiard system, see for example \cite{BT_anti}. We define two 1-forms $\nu = X^*\mu$,
$\hat\nu = X^*\hat\mu$ on $U\subset\mC^n$. In more detail,
$$
    \nu
  = \partial_2 L(\chi\!\circ\!\rho^{-1},\chi)\, d\chi, \quad
    \hat\nu
  = \partial_1 L(\chi\!\circ\!\rho^{-1},\chi)\, d\chi\!\circ\!\rho^{-1}.
$$
Equation (\ref{conj_Lagr}) is equivalent to
\begin{equation}
\label{nunu}
  \nu + \rho^*\hat\nu = 0.
\end{equation}

\begin{lem}
\label{lem:exact}
Equation (\ref{nunu}) implies
\begin{equation}
\label{<L>}
  \langle L(\chi,\chi\circ\rho) \rangle = L(0,0).
\end{equation}
\end{lem}

{\it Proof}. By (\ref{nunu})
$$
    0
  = \langle \nu + \rho^*\hat\nu \rangle
  = \langle \nu + \hat\nu \rangle
  = \langle dL(\chi\circ\rho^{-1},\chi) \rangle.
$$
It remains to use identity (\ref{<>=<>}). \qed

\subsection{Main equation}

Instead of (\ref{conj_Lagr}) or (\ref{nunu}) we consider an equivalent system $[\nu + \rho^*\nu] = 0$, (\ref{<L>}), or, in an explicit form\footnote
{At first glance it looks strange that for operators $\tau_\pm$ acting differently on $\chi$ and $f\circ\chi$ (the sign at the second term differs) the same notation is used. However below we will see that the function $\chi$ is odd in $z,\bar z$ while $f\circ\chi$ is even. Hence $\tau_\pm$ admit a universal formula $\tau_\pm = \id + \imath^*\circ\rho^{\pm1*}$, where $\imath$ is the central symmetry $(z,\bar z)\mapsto(-z,-\bar z)$.}
\begin{eqnarray}
\label{[explicit]}
&\displaystyle
   \Big[ \Big(
    \frac{\tau_-\chi + \tau_- f\circ\chi \, (\nabla f)\circ\chi}
       {L(\chi\circ\rho^{-1},\chi)}
  + \frac{\tau_+\chi + \tau_+ f\circ\chi \, (\nabla f)\circ\chi}
         {L(\chi,\chi\circ\rho)} \Big) d\chi \Big]
  = 0, & \\
\label{<L>+}
&   \Big\langle \big((\tau_-\chi)^2 + (\tau_- f\circ\chi)^2\big)^{1/2}
    \Big\rangle
  = L(0,0), & \\
\nonumber
&  \tau_\pm\chi = \chi - \chi\circ\rho^\pm, \quad
   \tau_\pm f\circ\chi = f\circ\chi + f\circ\chi\circ\rho^\pm.
&
\end{eqnarray}

For any $j = 1,\ldots,n$ consider the maps
$$
  \kappa,\kappa_j : B\to B, \quad
    \kappa_j (x_1,\ldots,x_n)
  = \big((-1)^{\delta_{j1}}x_1,\ldots,(-1)^{\delta_{jn}}x_n\big), \quad
    \kappa = \kappa_1\circ\ldots\circ\kappa_j.
$$
Any map $\kappa_j$ is the symmetry in the $j$-th coordinate hyperplane and $\kappa$ is the central symmetry w.r.t. the origin.

Given the frequency vector (\ref{lambda}) we search for a solution $(f,\chi)$ of (\ref{[explicit]}),(\ref{<L>+}) satisfying the following properties.
\medskip

{\bf (1)}. The function $f$ is real and (totally) even:
$$
  f\circ\kappa_j = f\quad\mbox{for any $j=1,\ldots,n$}.
$$
Its expansion in homogeneous forms $f^{(2k)}$ of degree $2k$ is
\begin{eqnarray*}
& f = f^{(0)} + f^{(2)} + f^{(4)} + \ldots, \qquad
  f^{(0)} = f_0 < 0, \quad
  f^{(2k)}(x) = \sum_{s\in\mZ_+^n,\,\|s\|=k} \FF_{2s} x^{2s}, & \\
& \FF_{2s}\in\mR, \quad
  x^{2s} = x_1^{2s_1} \ldots x_n^{2s_n}, \quad
  \|s\| = s_1 + \ldots + s_n.
&  
\end{eqnarray*}

{\bf (2)}. The conjugacy map $X$ is real:
\begin{equation}
\label{chi_real}
  \chi_-(z,\bar z) = \overline{\chi_-(z,\bar z)}, \quad
  \chi(z,\bar z) = \overline{\chi(z,\bar z)}.
\end{equation}

{\bf (3)}. The functions $\chi = (\chi_1,\ldots,\chi_n)$ are odd:
$$
  \chi_j\circ\kappa = - \chi_j, \qquad
  j=1,\ldots,n.
$$
Their expansions in homogeneous forms $\chi_j^{(k)}$ of degree $k$ is
$$
    \chi_j
  = \chi_j^{(1)} + \chi_j^{(3)} + \ldots, \qquad
   \chi^{(2k+1)}(z,\bar z)
  = \sum_{\|s'\|+\|s''\|=2k+1}^\infty \chi_{s's''} z^{s'} \bar z^{s''}.
$$

Note that two equations (\ref{chi_real}) are equivalent to each other and imply
$$
  \chi_{s's''} = \bar \chi_{s''s'}\quad
  \mbox{for any $s',s''\in\mZ^n_+$}.
$$  

{\bf (4)}. 
$\det\frac{\partial(\chi_-,\chi)}{\partial(z,\bar z)}\big|_{z=\bar z=0} \ne 0$.

\section{Symmetries}
\label{sec:symm}

\begin{dfn}
The rotation vector $\lambda$ is said to be nonresonant if the equation
$$
  \lambda_1^{k_1}\ldots \lambda_n^{k_n} = 1, \qquad
  k_1,\ldots,k_n\in\mZ
$$
holds only in the trivial situation: $k_1 = \ldots = k_n = 0$.
\end{dfn}

In this section we study symmetries of equation (\ref{[explicit]}),(\ref{<L>+}) with conditions {\bf (1)}--{\bf (4)}.
\medskip
 
{\bf (a)} (Gauge symmetry). System (\ref{[explicit]}), (\ref{<L>+}), 
{\bf (1)}--{\bf (4)} admits the following {\it gauge symmetry}. If $(f,\chi)$ is a solution, the pair $(f,\chi\circ s)$ is also a solution for any real $s : U\to U$ which commutes with $\rho$. Reality of the vector-function $s$ means that
\begin{equation}
\label{AAAA}
    s(z,\bar z)
  = \big( A_1(z,\bar z),\ldots, A_n(z,\bar z),
          \bar A_1(z,\bar z),\ldots, \bar A_n(z,\bar z)
    \big).
\end{equation}

\begin{lem}
Suppose that $\lambda$ is nonresonant.
Suppose that $s$, commuting with $\rho$, is real and can be expanded into a power series at zero. Then
\begin{equation}
\label{s}
     s(z,\bar z)
  = \big( z_1\psi_1,\ldots, z_n\psi_n,
          \bar z_1\bar\psi_1,\ldots,\bar z_n\bar\psi_n
    \big), \qquad
    \psi_j = \psi_j(z_1\bar z_1,\ldots,z_n\bar z_n).
\end{equation}
\end{lem}

{\it Proof}. We put
$$
  s = \sum s_{uv} z^u \bar z^v, \qquad
  s_{uv} = (A_{uv1},\ldots,A_{uvn},B_{uv1},\ldots,B_{uvn}), \quad
  u,v\in\mZ^n_+.
$$
Then the equation $s\circ\rho = \rho\circ s$ for any $j = 1,\ldots,n$ implies
$$
        (\lambda^{u-v} - \lambda_j) A_{uvj}
    =   0, \quad
        (\lambda^{u-v} - \bar\lambda_j) B_{uvj}
    =   0,  \qquad
        \lambda^{u-v}
    =   \lambda_1^{u_1-v_1} \ldots \lambda_n^{u_n-v_n}.
$$
By nonresonant condition we obtain
$$
  A_j(z,\bar z) = z_j\psi_j(z_1\bar z_1,\ldots,z_n\bar z_n), \quad
  B_j(z,\bar z) = \bar z_j\thet_j(z_1\bar z_1,\ldots,z_n\bar z_n).
$$
It remains to use (\ref{AAAA}). \qed
\medskip

{\bf (b)}  Consider in (\ref{<L>+}) the homogeneous form of degree two in
$z,\bar z$:
\begin{equation}
\label{deg1}
     \Big\langle
       \frac1{4f_0} \sum_j (\tau_- \chi_j^{(1)})^2 
      + \tau_- f^{(2)}\circ \chi^{(1)} \Big\rangle
   =  0. 
\end{equation}
By using the notation
$$
  \chi_j^{(1)} = \sum_{l=1}^n (c_{jl} z_l + \bar c_{jl} \bar z_l)
$$
after simple transformations we obtain:
$$
  \sum_{j,l} \big( 2 - \lambda_l - \lambda_l^{-1} + 8f_0 \FF_{e_j} \big)
              |c_{jl}|^2 z_l\bar z_l
  = 0,
$$
where $e_j\in\mZ_+^n$ is the $j$-th unit vector: its $j$-th components equals 1 while all others vanish.

This means that for any $l=1,\ldots,n$
$$
    \big( 2 - \lambda_l - \lambda_l^{-1} + 8f_0 \FF_{e_j} \big) |c_{jl}|^2
  = 0.
$$

Since $\lambda_j\ne\lambda_k^{\pm1}$ for $j\ne k$, we see that for any $l$ only one of coefficients $c_{jl}$ may be nonzero, the one, corresponding to $j$ such that $2 - \lambda_l - \lambda_l^{-1} + 8f_0\FF_j = 0$. By ${\bf (4)}$ for any $j$ such $j = j(l)$ exists. Without loss of generality we have:
$j(l) = l$. Hence,
\begin{equation}
\label{f_2}
    - 8f_0\FF_{e_j}
  = 2 - \lambda_j - \lambda_j^{-1}, \quad
    \chi_j^{(1)}
  = c_{jj} z_j + \bar c_{jj} \bar z_j, \quad
    c_{jj} \ne 0.
\end{equation}

Direct computation shows that by equations (\ref{f_2}) the homogeneous forms of degree 2 in (\ref{[explicit]}) vanish:
\begin{equation}
\label{deg2}
  \frac{(\tau_- + \tau_+)\chi_j^{(1)}
              + 4f_0\, \nabla_j f^{(2)}\circ\chi^{(1)}}
       {2|f_0|}
  \, d\chi_j^{(1)}
  = 0.
\end{equation}

We can assume that the coefficients $c_{jj}$ are real and positive. Indeed, by using the map $s$ (\ref{s})
$$
    s(z,\bar z)
  = (|c_{11}| c_{11}^{-1} z_1,\ldots,|c_{nn}| c_{nn}^{-1} z_n,
     |c_{11}| \bar c_{11}^{-1} \bar z_1,\ldots,
                                  |c_{nn}| \bar c_{nn}^{-1} \bar z_n)
$$
for the gauge transformation $(f,\chi)\mapsto (f,\chi\!\circ\! s)$, we obtain:
\begin{equation}
\label{chi^1=}
    (\chi_j\!\circ\! s)^{(1)}
  = |a_j| (z_j + \bar z_j), \qquad
    a_j = |c_{jj}|.
\end{equation}
\medskip

{\bf (c)} We define the maps (complex conjugacy of one coordinate)
$$
  \varkappa_1,\ldots,\varkappa_n:\mC^n\to\mC^n,\quad
  \varkappa_j(z) = w, \qquad
  w_j = \bar z_j, \quad
  w_k = z_k \;
  \mbox{ for any $k\ne j$}.
$$
Let $\rho_j : \mC^n\to\mC^n$ be the map
$$
  z \mapsto w = \rho_j(z), \qquad
  w_l = (\varkappa_j\lambda)_l z_l,
$$
where $(\varkappa_j\lambda)_l$ is the $l$-th coordinate of the vector $\varkappa_j\lambda$.
We have the identity
$$
  \rho\circ\varkappa_j = \varkappa_j\circ\rho_j.
$$

\begin{cor}
The pair $(f,\chi)$ is a solution of system (\ref{[explicit]}),(\ref{<L>+}), corresponding to the frequency vector $\lambda$ iff $(f,\chi\circ\varkappa_j)$ is a solution of system (\ref{[explicit]}),(\ref{<L>+}), corresponding to the frequency vector $\varkappa_j\lambda$.
\end{cor}

\begin{cor}
Without loss of generality it is possible to assume that
$\arg\lambda_j\in (0,\pi)$ for any $j = 1,\ldots,n$.
\end{cor}
\medskip

{\bf (d)} We define the maps $\TT_{jl}:\mC^n\to\mC^n$ which exchange
the coordinates $z_j$ and $z_l$ in any vector $z\in\mC^n$.
Let $\rho_{jl} : \mC^n\to\mC^n$ be the map
$$
  z \mapsto w = \rho_{jl}(z), \qquad
  w_k = (\TT\lambda)_k z_k, \qquad
  k = 1,\ldots,n.
$$
We have the identity
$$
  \rho\circ\TT_{jl} = \TT_{jl}\circ\rho_{jl}.
$$

\begin{cor}
The pair $(f,\chi)$ is a solution of system (\ref{[explicit]}),(\ref{<L>+}), corresponding to the frequency vector $\lambda$ iff $(f,\chi\circ\TT_{jl})$ is a solution of system (\ref{[explicit]}),(\ref{<L>+}), corresponding to the frequency vector $\TT_{jl}\lambda$.
\end{cor}

\section{Formal solution}
\label{sec:formal}

\begin{theo}
\label{theo:formal}
Suppose that the frequency vector $\lambda$ (\ref{lambda}) is nonresonant. Then for any $f_0<0$ and $a_1,\ldots,a_n > 0$ system (\ref{[explicit]}),(\ref{<L>+}) {\bf (1)}--{\bf (4)} has a formal solution given by power series
\begin{eqnarray}
\label{f2k}
&\displaystyle
  f = \sum_{k=0}^\infty f^{(2k)}, \qquad
  f^{(2k)}(x) = \sum_{\|s\|=k} f_{2s} x^{2s}, \quad
  \chi = \sum_{k=0}^\infty \chi^{(2k+1)}, & \\
\label{f2chi1}
&\displaystyle
    f^{(2)}
  = \sum_{j=1}^n \FF_{e_j} x_j^2, \quad
    \chi^{(1)}_j
  = a_j(z_j + \bar z_j), \qquad
    \FF_{e_j}
  = \frac{2 - \lambda_j - \lambda_j^{-1}}{-8f_0}. &
\end{eqnarray}
\end{theo}

{\it Proof of Theorem \ref{theo:formal}}. According to (\ref{f2k})--(\ref{f2chi1}) we are looking for $f$ even in $x$ and $\chi$ odd in $z,\bar z$. Hence even in $z,\bar z$ homogeneous forms in (\ref{[explicit]}),(\ref{<L>+}) will vanish. The forms of degree 2 in (\ref{[explicit]}) and (\ref{<L>+}) have been analyzed in item
{\bf (b)}, Section \ref{sec:symm}. In this way we obtain (\ref{f2chi1}).

Suppose that we have computed $f^{(2m)}$ and $\chi^{(2m-1)}$ for all $m<k$ by considering homogeneous forms in (\ref{[explicit]}),(\ref{<L>+}) of degrees $2,4,\ldots,2k-2$.

Taking in (\ref{[explicit]}),(\ref{<L>+}) the homogeneous form of degree
$2k$, we obtain:
\begin{eqnarray*}
&&\!\!\!\!\!\!\!\!\!\!\!\!\!\!\!\!\!\!\!\!
   \bigg[
     \frac{ (\tau_- + \tau_+) \chi_j^{(2k-1)}
           + 2f_0 (\nabla_j f^{(2k)}\circ\chi^{(1)}
                  + \nabla_j f^{(2)}\circ\chi^{(2k-1)}) }
          { 2|f_0| }
     \, d\chi_j^{(1)} \\
&&\qquad\qquad\quad     
  +\, \frac{ (\tau_- + \tau_+) \chi_j^{(1)}
           + 2f_0 \nabla_j f^{(2k)}\circ\chi^{(1)} }
          { 2|f_0| }
     \, d\chi_j^{(2k-1)}  \bigg] 
 \, = \, R_j^{(2k)},  \\
&&\!\!\!\!\!\!\!\!\!\!\!\!\!\!\!\!\!\!\!\!\!\!\!\!\!\!\!\!\!
   \Big\langle \frac 1{4f_0} \sum_{j=1}^n
            2\tau_-\chi_j^{(1)}\, \tau_-\chi_j^{(2k-1)}
          + \tau_- \sum_{j=1}^n 2\FF_{e_j} \chi_j^{(1)}\chi_j^{(2k-1)} 
          + \tau_- f^{(2k)}\circ\chi^{(1)}
   \Big\rangle
 \, = \, Q_j^{(2k)}, 
\end{eqnarray*}
where the forms $R_j^{(2k)}$ and the functions $Q_j^{(2k)}$ are polynomials w.r.t. coefficients of $f^{(2m)}$ and $\chi^{(2m-1)}$ with $m<k$,
$$
  R_j^{(2k)} = [ R_j^{(2k)} ], \quad
  Q_j^{(2k)} = \langle Q_j^{(2k)} \rangle.
$$
By (\ref{deg2}) the second fraction in the brackets $[\;]$ vanishes. Therefore
\begin{eqnarray}
\label{no1}
\!\!\!\!\!\!\!\!
   \Big[\Big(
      (\lambda_j + \lambda_j^{-1} - \rho^* - \rho^{-1*}) \chi_j^{(2k-1)}
           + 2f_0 \nabla_j f^{(2k)}\circ\chi^{(1)} \Big)
     \, d\chi_j^{(1)} \Big]
 &=& 2|f_0| R_j^{(2k)},  \\
\nonumber
\!\!\!\!\!\!\!\!
   \Big\langle \sum_{j=1}^n \tau_-\chi_j^{(1)}\, \tau_-\chi_j^{(2k-1)}
           + 8f_0 \sum_{j=1}^n \FF_{e_j} \chi_j^{(1)}\chi_j^{(2k-1)}
   \qquad\qquad\qquad && \\
\label{no2}           
           + 4f_0 \sum_{s\in\mZ_+^n,\,\|s\| = k} \FF_{2s} (\chi^{(1)})^{2s}
   \Big\rangle
 &=& 2f_0 Q_j^{(2k)},
\end{eqnarray}

Consider equation (\ref{no2}). The first term in the left-hand side equals
\begin{eqnarray*}
&&\Big\langle \sum_{j=1}^n
    a_j \Big( (1-\lambda_j^{-1}) z_j + (1-\lambda_j) \bar z_j \Big)
    \sum_{\|l'\| + \|l''\| = 2k-1}
     \chi_{j\,l'\,l''} (1-\lambda^{l''-l'}) z^{l'} \bar z^{l''} \Big\rangle \\
&=& \sum_{j=1}^n a_j(1-\lambda_j^{-1})(1-\lambda_j)
     \sum_{\|l\| = k} (\chi_{j\,l-e_j\,l} + \chi_{j\,l\,l-e_j}) z^l\bar z^l.
\end{eqnarray*}

The second term equals
\begin{eqnarray*}
&&\Big\langle 8f_0\sum_{j=1}^n
    \FF_{e_j} a_j (z_j + \bar z_j)
    \sum_{\|l'\| + \|l''\| = 2k-1}
       \chi_{j\,l'\,l''} z^{l'} \bar z^{l''} \Big\rangle \\
&=& 8f_0\sum_{j=1}^n
    \frac{2-\lambda_j-\lambda_j^{-1}}{-8f_0} a_j
    \sum_{\|l'\| + \|l''\| = 2k-1}
       (\chi_{j\,l-e_j\,l} + \chi_{j\,l\,l-e_j}) z^l \bar z^l.
\end{eqnarray*}
Hence the first two terms in the left-hand side of (\ref{no2}) cancel.

The third term equals
$$
  4f_0 \sum_{s\in\mZ_+^n,\,\|s\| = k} \FF_{2s} a^{2s}
            C_{2s_1}^{s_1} \ldots C_{2s_n}^{s_n}
            (z_1\bar z_1)^{s_1} \ldots (z_n\bar z_n)^{s_n}.
$$
Hence the coefficients $\FF_{2s}$ are uniquely computed from (\ref{no2}).

Now turn to equation (\ref{no1}). The first term in the left-hand side equals
$$
  \sum_{\|l'\|+\|l''\|=2k-1,\,l'-l''\ne\pm e_j}
     (\lambda_j + \lambda_j^{-1} - \lambda^{l'-l''} - \lambda^{l''-l'})
       \chi_{j\,l'\,l''} z^{l'} \bar z^{l''}.
$$
The condition $l'-l''\ne\pm e_j$ appears as a result of application of the operation $[\;]$. It implies that the coefficients
$\lambda_j + \lambda_j^{-1} - \lambda^{l'-l''} - \lambda^{l''-l'}$ do not vanish.

The second term in the left-hand side has been already computed from equation (\ref{no2}). Hence the coefficients $\chi_{j\,l'\,l''}$, $l'-l''\ne\pm e_j$ are computed uniquely from (\ref{no1}).

The coefficients $\chi_{j\,l-e_j\,l}$ and
$\chi_{j\,l\,l-e_j} = \bar\chi_{j\,l-e_j\,l}$ can be chosen arbitrarily due to the gauge symmetry.
 \qed

\section{Other symmetries}
\label{sec:other}

When the existence theorem (Theorem \ref{theo:formal}) is proven, we can return to equation (\ref{conj_Lagr}) which is equivalent to system (\ref{[explicit]}),(\ref{<L>+}) but looks somewhat simpler. Now we can ignore overdeterminacy of (\ref{conj_Lagr}) w.r.t. $f$. In an explicit form (\ref{conj_Lagr}) looks as follows:
\begin{equation}
\label{explicit}
  \frac{\tau_-\chi + \tau_- f\circ\chi \, (\nabla f)\circ\chi}
       {L(\chi\circ\rho^{-1},\chi)}
  + \frac{\tau_+\chi + \tau_+ f\circ\chi \, (\nabla f)\circ\chi}
         {L(\chi,\chi\circ\rho)}
  = 0.
\end{equation}

Direct computation shows that, analogously to (\ref{no1}), $f^{(2k)}$ and $\chi^{(2k-1)}$ satisfy the equations
\begin{equation}
\label{no1+}
\!\!\!\!\!\!\!\!\!
       (\lambda_j + \lambda_j^{-1} - \rho^* - \rho^{-1*}) \chi_j^{(2k-1)}
           + 2f_0 \nabla_j f^{(2k)}\circ\chi^{(1)}
  =  2|f_0| P_j^{(2k)},
\end{equation}
where the functions $P_j^{(2k)}$ are polynomials w.r.t. coefficients of $f^{(2m)}$ and $\chi^{(2m-1)}$ with $m<k$.

\begin{prop}
(1) If $\chi_j$ are chosen satisfying (\ref{f2chi1}) then
$\chi_j$ are odd in $z_j,\bar z_j$ and even in $z_l,\bar z_l$ for all $l\ne j$
\begin{equation}
\label{chikappa}
  \chi_j = -\chi_j\circ\kappa_j = \chi_j\circ\kappa_l\;
   \mbox{ for all $l\ne j$}.
\end{equation}

(2) If the coefficients $\chi_{j\,l-e_j,\,l} = \bar\chi_{j\,l,\,l-e_j}$ (arbitrary due to the gauge symmetry) are chosen real, then all Taylor coefficients of $\chi$ are real.
\end{prop}

{\it Proof}. (1) It is sufficient to note that if the homogeneous forms $\chi^{(2m-1)}$ with $m<k$ satisfy (\ref{chikappa}), the functions $P_j^{(2k)}$ also satisfy (\ref{chikappa}):
$$
   P_j^{(2k)} = -P_j^{(2k)}\circ\kappa_j = P_j^{(2k)}\circ\kappa_l\;
   \mbox{ for all $l\ne j$}.
$$

(2) If the homogeneous forms $\chi^{(2m-1)}$ with $m<k$ satisfy (\ref{chi_real}) and coefficients of the forms $\chi^{(2m-1)}$ are real, the functions $P_j^{(2k)}$ are also real polynomials in $z$ and $\bar z$ with real coefficients. \qed
\medskip

The case $\lambda_1 = \ldots = \lambda_n = -1$ is resonant and one should not expect a positive answer to question ${\bf\hat Q}$. Although in this case $\chi$ can not be found, the function $f$ can be computed explicitly. Indeed, in this case
$$
  \tau_\pm \chi_j = 2\chi_j, \quad
  \tau_\pm f\circ\chi = 2f\circ\chi, \quad
    \hat L(\chi\circ\rho^{-1},\chi)
  = \hat L(\chi,\chi\circ\rho)
  = 2\sqrt{h}, \qquad
    h
  = f^2\circ\chi + \chi^2.
$$

In this case in equations (\ref{explicit}) we can regard $\chi_j$ as independent variables. These equations are equivalent to
$\partial_{\chi_j} h = 0$. Condition {\bf (1)} implies that $h = f_0^2$, therefore $f(x) = -\sqrt{f_0^2 - x^2}$.

This formal computation shows that, in a certain sense, any sphere with the center at the origin is a limit solution of our problem when all $\lambda_j$ tend to $-1$. Numeric computations confirm this statement.

\section{Another version of the main equation}
\label{sec:another}

In this section we present another version of the equations from which the functions $f$ and $\chi$ can be computed. For $n=1$ these equations are the same as in \cite{Tre_PhysD,Tre_Proc}. We used these equations for $n=2$ in numeric computations since they contain the unknown functions polynomially.

\subsection{Dynamical equations}

Let $\ba,\bb,\bc\in E$ be three successive points of a billiard trajectory:
\begin{equation}
\label{abc}
  \ba = \Big(\begin{array}{c} a \\ -f(a) \end{array} \Big), \quad
  \bb = \Big(\begin{array}{c} b \\ f(b) \end{array} \Big), \quad
  \bc = \Big(\begin{array}{c} c \\ -f(c) \end{array} \Big).
\end{equation}
The vector $\bn$ normal to $S_+$ at the point $b$ is
$$
  \bn = \Big( \begin{array}{c} s\\ -1\end{array} \Big), \qquad
  s_j = \partial_{x_j} f(b), \quad j = 1,\ldots,n.
$$
The law of elastic reflection implies
$$
                \ba-\bb
  \;\parallel\; \bc-\bb - 2\bn\frac{\langle \bn,\bc-\bb\rangle}{\bn^2}
$$
i.e., the vector $\ba-\bb$ is collinear to the vector symmetric to $\bc-\bb$ w.r.t. a plane orthogonal to $\bn$. In more detail,
$$
               \ba - \bb
  \;\parallel\; A (\bc - \bb),   \qquad
    A
  = \bn^2 I_{n+1} - 2 \bn \bn^T = C^{-1} \bI C,
$$
where
$$
    \bI
  = \operatorname{diag}(1,\ldots,1,-1), \quad
    C
  = \left(\begin{array}{ccccc}  0 &   0   &\ldots &   1   & s_n  \\
                           \ldots &\ldots &\ldots &\ldots &\ldots  \\
                                0 &   1   &\ldots &   0   & s_2  \\
                                1 &   0   &\ldots &   0   & s_1  \\
                              s_1 &  s_2  &\ldots &  s_n  & -1
      \end{array}\right) .
$$
Therefore we obtain:
\begin{equation}
\label{dynamic}
                C (\bb - \ba)
 \;\parallel\;  \bI C (\bb - \bc).
\end{equation}
We refer to the map $(a,b) \mapsto\hat\beta(a,b) = (b,c)$ as the billiard map.

\subsection{Conjugacy equation}
\label{sec:conj_old}

Question $\bf{\hat Q}$ is equivalent to the following one. Does the following conjugacy equation
\begin{equation}
\label{conj0}
  \hat\beta\circ X = X\circ\rho, \qquad
  \hat\beta = \hat\beta(f)
\end{equation}
have a solution $(f,X)$, where $f$ satisfies (\ref{f}) and
$X = (\chi_-,\chi) : U\to B^2$ is a diffeomorphism of a neighborhood $U\subset\mR^{2n}$ to its image.

Equations (\ref{conj0}), (\ref{chichi}), and (\ref{dynamic}) imply
\begin{equation}
\label{conj_n}
               C\tau_- \Big(\begin{array}{c}\chi\\ f\circ\chi\end{array}\Big)
 \;\parallel\; \bI C\tau_+
                   \Big(\begin{array}{c}\chi\\ f\circ\chi\end{array}\Big).
\end{equation}

For $n=1$ equation (\ref{conj_n}) takes the form
$$
         \Big(\begin{array}{c} \tau_-\chi + s\,\tau_- f\circ\chi \\
                               s\,\tau_-\chi - \tau_- f\circ\chi
              \end{array}\Big)
 \;\parallel\;
         \Big(\begin{array}{c} \tau_+\chi + s\,\tau_+ f\circ\chi \\
                              -s\,\tau_+\chi + \tau_+ f\circ\chi
              \end{array}\Big),
$$
which implies
\begin{eqnarray}
\nonumber
&&\!\!\!\!\!\!
   2f'\circ\chi \big( \tau_- f\circ\chi\, \tau_+ f\circ\chi
                      - \tau_-\chi\, \tau_+\chi \big) \\
\label{conj_eq_old}
&&  + \big( 1-f^{\prime\, 2}\circ\chi \big)
                \big( \tau_- \chi\, \tau_+ f\circ\chi
                       + \tau_+\chi\, \tau_- f\circ\chi \big)
  = 0.
\end{eqnarray}
This equation is obtained in \cite{Tre_PhysD}. Its analysis is contained in \cite{Tre_PhysD,Tre_Proc}.

If $n=2$, (\ref{conj_n}) takes the form
$$
         \Bigg(\begin{array}{c} \tau_-\chi_2 + s_2\,\tau_- f\circ\chi \\
                               \tau_-\chi_1 + s_1\,\tau_- f\circ\chi \\
                  s_1\,\tau_-\chi_1 + s_2\,\tau_-\chi_2 - \tau_- f\circ\chi
              \end{array}\Bigg)
 \;\parallel\;
         \Bigg(\begin{array}{c} \tau_+\chi_2 + s_2\,\tau_+ f\circ\chi \\
                               \tau_+\chi_1 + s_1\,\tau_+ f\circ\chi \\
                 - s_1\,\tau_+\chi_1 - s_2\,\tau_+\chi_2 + \tau_+ f\circ\chi
              \end{array}\Bigg),
$$
which implies 2 equations:
\begin{eqnarray}
\nonumber
&&\!\!\!\!
    \, 2\partial_1 f\circ\chi
     \big( \tau_- f\circ\chi\,\tau_+ f\circ\chi
          - \tau_- \chi_1\, \tau_+ \chi_1 \big) \\
\nonumber
&&\!\!\!\!\!\!\!\!\!\!\!
  + \, \big( 1 - (\partial_1 f\circ\chi)^2 \big)
      \big( \tau_-\chi_1\,\tau_+ f\circ\chi
           + \tau_+\chi_1\,\tau_- f\circ\chi \big) \\
\nonumber
&&\!\!\!\!\!\!\!\!\!\!\!
  - \, \partial_2 f\circ\chi
     \big( \tau_+\chi_2\,\tau_-\chi_1 + \tau_-\chi_2\, \tau_+\chi_1 \big) \\
\label{conj_2}
&&\!\!\!\!\!\!\!\!\!\!\!
  - \, \partial_2 f\circ\chi\,\partial_1 f\circ\chi
    \big( \tau_-\chi_2\, \tau_+ f\circ\chi
         + \tau_+\chi_2\, \tau_- f\circ\chi \big)
 \, = \, 0, \\[1.5mm]
\nonumber
&&\!\!\!\!
    \, 2\partial_2 f\circ\chi
     \big( \tau_- f\circ\chi\,\tau_+ f\circ\chi
          - \tau_- \chi_2\, \tau_+ \chi_2 \big) \\
\nonumber
&&\!\!\!\!\!\!\!\!\!\!\!
  + \, \big( 1 - (\partial_2 f\circ\chi)^2 \big)
      \big( \tau_-\chi_2\,\tau_+ f\circ\chi
           + \tau_+\chi_2\,\tau_- f\circ\chi \big) \\
\nonumber
&&\!\!\!\!\!\!\!\!\!\!\!
  - \, \partial_1 f\circ\chi
     \big( \tau_+\chi_1\,\tau_-\chi_2 + \tau_-\chi_1\, \tau_+\chi_2 \big) \\
\label{conj_1}
&&\!\!\!\!\!\!\!\!\!\!\!
  - \, \partial_1 f\circ\chi\,\partial_2 f\circ\chi
    \big( \tau_-\chi_1\, \tau_+ f\circ\chi
         + \tau_+\chi_1\, \tau_- f\circ\chi \big)
 \, = \, 0.
\end{eqnarray}

Our numeric results in the case $n=2$ presented in Section \ref{sec:intro} are based on the analysis of system (\ref{conj_2}),(\ref{conj_1}).

\end{document}